May 25, 2022

# Stirling numbers and inverse factorial series

Khristo N Boyadzhiev

Department of Mathematics
Ohio Northern University, Ada, Ohio 45810, USA
k-boyadzhiev@onu.edu

**Abstract**

We study inverse factorial series and their relation to Stirling numbers of the first kind. We prove a special representation of the polylogarithm function in terms of series with such numbers. Using various identities for Stirling numbers of the first kind we construct a number of expansions of functions in terms of inverse factorial series where the coefficients are special numbers. These results are used to prove/reprove the asymptotic expansion of some classical functions. We also prove a binomial formula involving inverse factorials.

## 1. Introduction

The unsigned Stirling numbers of the first kind $\begin{bmatrix} n \\ k \end{bmatrix}$ first appeared in James Stirling's book

*Methodus Differentialis* (published in Latin, in 1730). See [26] for a translation and comments. These numbers have numerous applications in combinatorics and analysis ([1, 3, 6, 9, 10, 12, 24]).

In this paper we present various results involving the unsigned Stirling numbers of the first kind and inverse factorial series.

These numbers are usually defined by their ordinary generating function

$$(1) \qquad x(x+1)...(x+n-1) = \sum_{k=0}^{n} \begin{bmatrix} n \\ k \end{bmatrix} x^k$$

or by their exponential generating function

$$(2) \quad \frac{[-\ln(1-x)]^k}{k!} = \sum_{n=k}^{\infty} \begin{bmatrix} n \\ k \end{bmatrix} \frac{x^n}{n!}, \quad |x|<1, k \geq 0$$

(see [9, 10]). However, there is another, less known generating function coming from the original works of James Stirling

$$(3) \quad \frac{1}{x^{k+1}} = \sum_{n=k}^{\infty} \begin{bmatrix} n \\ k \end{bmatrix} \frac{1}{x(x+1)\ldots(x+n)}$$

in terms of inverse factorials (see [3] and [26, p.171]). Series of this form are called inverse factorial series. As shown by Milne-Thomson in Chapter 10 of his book [16], for appropriate functions $f(z)$ defined on a half-plane $\operatorname{Re}(z) > \delta > 0$ expansions of the form

$$f(z) = \sum_{n=0}^{\infty} \frac{a_n}{z(z+1)\ldots(z+n)}$$

are convergent and the coefficients are uniquely determined. See also Whittaker and Watson [30], pp. 142-144. Inverse factorial series were thoroughly studied by Norlund in [18, 19, 20], Obrechkoff [22], and Watson [27], and in more recent times by Wasow in [28, Chapter 11], Weniger [29], and Karp with Prilepkina [13]. Instead of listing general conditions for convergence, in many cases it is more practical to prove convergence by using the ratio test or the Raabe-Duhamel test.

Formula (3) brings immediately to an important result. Setting $z = m$, an integer, summing for $m = 1, 2, \ldots$ and changing the order of summation we come to the representation

$$\zeta(k+1) = \sum_{n=k}^{\infty} \begin{bmatrix} n \\ k \end{bmatrix} \sum_{m=1}^{\infty} \frac{1}{m(m+1)\ldots(m+n)}$$

where $\zeta(s)$ is Riemann's zeta function

$$\zeta(s) = \sum_{n=1}^{\infty} \frac{1}{n^s}, \quad \operatorname{Re}(s) > 1.$$

It is easy to compute that

$$\sum_{m=1}^{\infty} \frac{1}{m(m+1)\ldots(m+n)} = \frac{1}{n!n}$$

(see, for example, [4]) and we obtain an important representation of the zeta function

$$(4) \quad \zeta(k+1) = \sum_{n=k}^{\infty} \begin{bmatrix} n \\ k \end{bmatrix} \frac{1}{n!n}.$$

This result is not new, it was known to Charles Jordan at least in 1939 (see equation (6) on p. 166 and also equation (11) on p. 194 in his book [12]). Jordan even proved a more general formula ([12], p. 343); for another proof see [4]. Informative comments on equation (4) can be found in Adamchik's paper [1]. See also problems 46 and 47 in Chapter 1 of [24].

In this paper we extend the representation (4) to the Hurwitz zeta function

$$\zeta(s,a) = \sum_{n=0}^{\infty} \frac{1}{(n+a)^s}, \ \operatorname{Re}(s) > 1, \ a > 0$$

(see Proposition 1 below). The next natural step is to replace the Riemann zeta function by the polylogarithm

$$\operatorname{Li}_s(x) = \sum_{n=1}^{\infty} \frac{x^n}{n^s} \ (|x| \leq 1).$$

This is our main result (Theorem 1 in Section 2). It represents an important extension since $\operatorname{Li}_s(1) = \zeta(s)$.

In Section 3 we list a number of examples of functions represented by inverse factorial series. In Section 4 with the help of Proposition 6 and using Stirling number identities we generate new inverse factorial representations. Then in Section 5 again with the help of Proposition 6 we obtain several important asymptotic expansions.

Finally, in Section 6 we prove a formula relating a binomial sum to inverse factorial series with unsigned Stirling numbers of the first kind.

## 2. Main results

**Proposition 1.** *For the Hurwitz zeta function* $\zeta(s,a)$ *we have the representation*

(5) $\quad \zeta(k+1,a) = \Gamma(a) \sum_{n=k}^{\infty} \begin{bmatrix} n \\ k \end{bmatrix} \frac{1}{n\Gamma(n+a)} \quad (k \geq 1, a > 0)$

*which turns into* (4) *when* $a = 1$.

Proof. With $x = m + a, \ m = 0, 1, 2, \ldots$ we write (3) in the form

$$\frac{1}{(m+a)^{k+1}} = \sum_{n=k}^{\infty} \begin{bmatrix} n \\ k \end{bmatrix} \frac{1}{(m+a)(m+a+1)\ldots(m+a+n)}$$

then summing for $m = 0, 1, 2, \ldots$ we come to

$$\zeta(k+1,a) = \sum_{n=k}^{\infty} \begin{bmatrix} n \\ k \end{bmatrix} \left\{ \sum_{m=0}^{\infty} \frac{1}{(m+a)(m+a+1)\ldots(m+a+n)} \right\}.$$

Evaluating the sum in the braces gives

$$\sum_{m=0}^{\infty}\frac{1}{(m+a)(m+a+1)\ldots(m+a+n)}=\frac{1}{n\Gamma(n+a)}$$

and (5) is proved.

Next we recall the harmonic numbers

$$H_n=1+\frac{1}{2}+\ldots+\frac{1}{n},\ H_0=0$$

which appear in the next proposition.

**Proposition 2**. *For any* $n\geq 0$

(6) $$\sum_{p=1}^{\infty}\frac{x^{p+n}}{p(p+1)\ldots(p+n)}=P_n(x)+\frac{(-1)^{n-1}}{n!}(1-x)^n\ln(1-x)$$

*where* $P_n(x)$ *is the polynomial*

(7) $$P_n(x)=\frac{1}{n!}\left(H_n(x-1)^n+\sum_{j=0}^{n-1}\binom{n}{j}\frac{(x-1)^j}{n-j}\right)=\frac{1}{n!}\left(H_n(x-1)^n+\sum_{k=1}^{n}\binom{n}{k}\frac{(x-1)^{n-k}}{k}\right)$$

*with* $P_n(1)=\dfrac{1}{n!n}$.

The presence of the harmonic numbers $H_n$ here is very interesting!

Proof. Let us set

$$f_n(x)=\sum_{p=1}^{\infty}\frac{x^{p+n}}{p(p+1)\ldots(p+n)}$$

for $n=0,1,\ldots$ . Differentiating this function $n$ times we find

$$\left(\frac{d}{dx}\right)^n f_n(x)=\sum_{p=1}^{\infty}\frac{x^p}{p}=-\ln(1-x)$$

We see that $f_n(x)$ is the $n$-th antiderivative of $-\ln(1-x)$ satisfying the conditions (8) below. We will find a closed form evaluation for the series $f_n(x)$.

Clearly $(d/dx)f_m(x)=f_{m-1}(x)$. When $n=0$,

$$f_0(x)=\sum_{p=1}^{\infty}\frac{x^p}{p}=-\ln(1-x)$$

and $f_1, f_2, ... f_n$ are obtained from here by repeated integration. The constant of integration is taken zero each time so that according to (4)

(8) $\quad f_1(1) = \dfrac{1}{1!1}, \; f_2(1) = \dfrac{1}{2!2}, \; f_3(1) = \dfrac{1}{3!3}, ..., f_n(1) = \dfrac{1}{n!n}$ .

The process starts this way

$$f_1(x) = \int -\ln(1-x)dx = \int \ln(1-x)d(1-x) = (1-x)\ln(1-x) + x + C$$

$$= (1-x)\ln(1-x) + x$$

then

$$f_2(x) = \dfrac{3}{4}x^2 - \dfrac{x}{2} - \dfrac{1}{2}(1-x)^2 \ln(1-x)$$

etc. We evaluate each time at $x = 1$ by taking limits,

$$\lim_{x \to 1}(1-x)^n \ln(1-x) = 0, \quad n = 1, 2, ... \; .$$

Integrating repeatedly by parts we come to the representation

(9) $\quad f_n(x) = \displaystyle\sum_{p=1}^{\infty} \dfrac{x^{p+n}}{p(p+1)...(p+n)} = P_n(x) + \dfrac{(-1)^{n-1}}{n!}(1-x)^n \ln(1-x)$

where $P_n(x)$ is a polynomial of degree $n$. For convenience we write it in the form

$$P_n(x) = a_n(x-1)^n + a_{n-1}(x-1)^{n-1} + ... + a_1(x-1) + a_0 .$$

Obviously, $a_0 = f_n(1) = \dfrac{1}{n!n}$. To compute the other coefficients we repeatedly differentiate in (9) and set $x \to 1$. The coefficient $a_n$ requires a longer (but elementary) computation. Namely, differentiating equation (9) $n$ times and setting $x \to 1$ brings to

$$0 = n!a_n + \sum_{k=1}^{n} \binom{n}{k} \dfrac{(-1)^k}{k}$$

and now using the well-known identity

$$\sum_{k=1}^{n} \binom{n}{k} \dfrac{(-1)^{k-1}}{k} = H_n .$$

we find that $n!a_n = H_n$ as desired. Thus we come to the representation (7) which completes the proof of the proposition.

Setting $x = -1$ in (7) we compute that

$$P_n(-1) = \frac{(-1)^n 2^n}{n!} \sum_{k=1}^{n} \frac{1}{2^k k}$$

and we come to the interesting result:

**Corollary 3.** *For every* $n \geq 0$ *we have*

(10) $$\sum_{p=1}^{\infty} \frac{(-1)^p}{p(p+1)...(p+n)} = \frac{2^n}{n!} \left( \sum_{k=1}^{n} \frac{1}{2^k k} - \ln 2 \right)$$

Now we come to our main result, a theorem that extends the representation (4) to polylogarithms.

**Theorem 4.** *For every integer* $k \geq 0$, *and* $|x| < 1$, *we have*

(11) $$\operatorname{Li}_{k+1}(x) = \sum_{n=k}^{\infty} \begin{bmatrix} n \\ k \end{bmatrix} \frac{1}{x^n} \left\{ P_n(x) + \frac{(-1)^{n-1}}{n!} (1-x)^n \ln(1-x) \right\}$$

where $P_n(x)$ is the polynomial from (7).

Note that when we set $x \to 1$ in (11) the equation turns into (4), as $P_n(1) = \frac{1}{n!n}$.

*Proof.* We write (3) in the form

$$\frac{1}{p^{k+1}} = \sum_{n=k}^{\infty} \begin{bmatrix} n \\ k \end{bmatrix} \frac{1}{p(p+1)...(p+n)}$$

then multiply both sides by $x^p$ and sum for $p$ from 1 to infinity. After changing the order of summation we write

$$\operatorname{Li}_{k+1}(x) = \sum_{n=k}^{\infty} \begin{bmatrix} n \\ k \end{bmatrix} \frac{1}{x^n} \left\{ \sum_{p=1}^{\infty} \frac{x^{p+n}}{p(p+1)...(p+n)} \right\}$$

and the result follows now from Proposition 2. The proof is completed.

**Remark**. A mentioned above, equation (9) represents the $n$th antiderivative of $-\ln(1-x)$ with the conditions (8). The functions $f_2(x)$ and $f_3(x)$ can be found in the table [21] of Prudnikov et al.; $f_2$ is entry 5.2.6 (3) and $f_3$ is entry 5.2.6(12). They are listed also in Hansen's table [11], $f_2$ is entry 5.7.40 and $f_3$ is entry 5.6.29. The table [11] has some general formulas similar to (9) like 10.4.15 and 10.4.17. Hansen refers to the works of Schwatt (see [23], pp.191-192). The functions $f_n$ were studied also by Mathar [14] and a similar study can be found in [15].

We finish this section with a related result. Consider the Dirichlet series

$$H(s) = \sum_{n=1}^{\infty} \frac{H_n}{n^s} \quad (\operatorname{Re} s > 1).$$

It is reasonable to ask the same question: can we express $H(k+1)$ in terms of $\begin{bmatrix} n \\ k \end{bmatrix}$? The affirmative answer is given in the following proposition.

**Proposition 5.** *For every* $k \geq 1$

(12) $$\sum_{p=1}^{\infty} \frac{H_p}{p^{k+1}} = \sum_{n=k}^{\infty} \begin{bmatrix} n \\ k \end{bmatrix} \frac{\psi'(n)}{n!}$$

*and also*

$$\left(\frac{k+3}{2}\right)\zeta(k+2) - \frac{1}{2}\sum_{j=1}^{k-1}\zeta(j+1)\zeta(k-j+1) = \sum_{n=k}^{\infty} \begin{bmatrix} n \\ k \end{bmatrix} \frac{\psi'(n)}{n!}.$$

*where* $\psi(z) = \Gamma'(z)/\Gamma(z)$ *is the digamma function.*

Proof. From (3) we have

$$\frac{H_p}{p^{k+1}} = \sum_{n=k}^{\infty} \begin{bmatrix} n \\ k \end{bmatrix} \frac{H_p}{p(p+1)...(p+n)}$$

and then we sum for $p = 1, 2, ...$. Interchanging the summations on the right hand side we write

$$\sum_{p=1}^{\infty} \frac{H_p}{p^{k+1}} = \sum_{n=k}^{\infty} \begin{bmatrix} n \\ k \end{bmatrix} \left\{ \sum_{p=1}^{\infty} \frac{H_p}{p(p+1)...(p+n)} \right\}.$$

It is known that for $n \geq 1$

$$\sum_{p=1}^{\infty} \frac{H_p}{p(p+1)...(p+n)} = \frac{1}{n!}\sum_{m=0}^{\infty} \frac{1}{(n+m)^2} = \frac{\psi'(n)}{n!}$$

(see, for instance, [4]). The Euler sums on the left hand side in (12) have the closed form evaluation [1]

$$\sum_{p=1}^{\infty} \frac{H_p}{p^{k+1}} = \left(\frac{k+3}{2}\right)\zeta(k+2) - \frac{1}{2}\sum_{j=1}^{k-1}\zeta(j+1)\zeta(k-j+1)$$

and completes the proof.

## 3. Examples of inverse factorial series

First we list several known series and then we add some new.

**Example 1.** Using the properties

$$\begin{bmatrix} n+1 \\ 1 \end{bmatrix} = n!, \quad \begin{bmatrix} n+1 \\ 2 \end{bmatrix} = n! H_n$$

we find from (3) with $k=1$ and $k=2$ the following two expansions

(13) $\quad \dfrac{1}{z-1} = \sum_{n=0}^{\infty} \dfrac{n!}{z(z+1)\ldots(z+n)}.$

$$\dfrac{1}{(z-1)^2} = \sum_{n=0}^{\infty} \dfrac{n! H_n}{z(z+1)\ldots(z+n)}$$

by making in (3) the substitutions $n \to n+1,\ z \to z-1$.

Both series are convergent for $\operatorname{Re}(z) > 1$. This convergence can be verified by the Raabe-Duhamel test. For instance, if $a_n$ is the general term of the series in (13) we have for real $z > 1$,

$$\lim n\left(\dfrac{a_n}{a_{n+1}} - 1\right) = z > 1.$$

The Raabe-Duhamel test confirms also the convergence in the following several expansions.

Let $w$ be a complex number. A series extending (13) is

(14) $\quad \dfrac{1}{z-w} = \sum_{n=0}^{\infty} \dfrac{w(w+1)\ldots(w+n-1)}{z(z+1)\ldots(z+n)}$

convergent for $\operatorname{Re}(z) > \operatorname{Re}(w)$ (see [19, p. 222]).

**Example 2.** Another example in the same spirit is the series

(15) $\quad \dfrac{1}{(z-1)(z-2)} = \sum_{n=0}^{\infty} \dfrac{n! n}{z(z+1)\ldots(z+n)}, \quad \operatorname{Re}(z) > 2$

(cf. entry 6.6.14 in Hansen's table [11]).

**Example 3.** Next (see 6.6.39 in [11]) we have

(16) $\quad \psi'(z) = \sum_{n=0}^{\infty} \dfrac{1}{(z+n)^2} = \sum_{n=0}^{\infty} \dfrac{n!}{n+1} \dfrac{1}{z(z+1)\ldots(z+n)}, \quad \operatorname{Re}(z) > 0$

where $\psi(z) = \Gamma'(z)/\Gamma(z)$ is the digamma function.

**Example 4**. Consider the beta function

$$\beta(z) = \int_0^1 \frac{t^{z-1}}{t+1} dt = \int_0^\infty \frac{e^{-zt}}{e^{-t}+1} dt = \sum_{n=0}^\infty \frac{(-1)^n}{n+z}, \quad \text{Re}(z) > 0.$$

This function was studied by Niels Nielsen [17] and is often called Nielsen's beta function. For some properties and important integrals related to $\beta(z)$ see [8].

Integration by parts gives

$$\beta(z) = \frac{1}{2z} + \frac{1}{z}\int_0^1 \frac{t^z}{(t+1)^2} dt$$

and repeating this again and again we come to the expansion

(17) $$\beta(z) = \sum_{n=0}^\infty \frac{n!}{2^{n+1}} \frac{1}{z(z+1)\ldots(z+n)}.$$

This representation appears on p. 81 in Nielsen's book [17]. It appears also on p. 352 in Norlund's paper [18].

**Example 5**. Consider the lower incomplete gamma function

$$\gamma(z, x) = \int_0^x t^{z-1} e^{-t} dt, \quad x > 0, \text{Re}(z) > 0.$$

We have (Temme [25, p. 279])

$$\gamma(z, x) = x^z e^{-x} \sum_{n=0}^\infty \frac{x^n}{z(z+1)\ldots(z+n)}.$$

Convergence follows from the ratio test.

**Example 6**. Jacques Binet proved the expansion for the log-gamma function

$$\ln \Gamma(z) = \left(z - \frac{1}{2}\right) \ln z - z + \ln\sqrt{2\pi} + \sum_{n=1}^\infty \frac{a_n}{(z+1)(z+2)\ldots(z+n)}$$

where

$$a_n = \frac{1}{n}\int_0^1 \left(t - \frac{1}{2}\right) t(t+1)\ldots(t+n-1) dt.$$

We can write this in the form (with $a_0 = 0$)

$$\frac{1}{z} \ln \frac{\Gamma(z)}{\sqrt{2\pi}} - \left(1 - \frac{1}{2z}\right) \ln z + 1 = \sum_{n=0}^\infty \frac{a_n}{z(z+1)(z+2)\ldots(z+n)}$$

(see Whittaker and Watson [30, p. 253]).

Inverse factorial series appeared also in the works of Ramanujan, as discussed by Berndt in [2]. More details will be given in Section 5.

## 4. Generating new inverse factorial representations

The following proposition makes it possible to generate more inverse factorial representations.

**Proposition 6**. *Let $a_0, a_1, ..., a_n, ...$ be a sequence. Then we formally have*

$$(18) \quad \sum_{k=0}^{\infty} \frac{a_k}{z^{k+1}} = \sum_{n=0}^{\infty} \frac{1}{z(z+1)...(z+n)} \left\{ \sum_{k=0}^{n} \begin{bmatrix} n \\ k \end{bmatrix} a_k \right\}.$$

For the proof we multiply both sides in (3) by $a_k$, sum for $k$ from zero to infinity and change the order of summation. We do not discuss convergence in the general case. Convergence can be checked easily in all particular examples below.

Using the Stirling sequence transformation

$$b_n = \sum_{k=0}^{n} \begin{bmatrix} n \\ k \end{bmatrix} a_k$$

we can generate various inverse factorial representations. A short table of Stirling transform identities can be found in the Appendix in [5]. The entries in this table are written in terms of the (signed) Stirling numbers of the first kind $s(n,k)$ for which

$$(19) \quad \begin{bmatrix} n \\ k \end{bmatrix} = (-1)^{n-k} s(n,k)$$

and in terms of $s(n,k)$ formula (18) can be written in the form

$$(20) \quad \sum_{k=0}^{\infty} \frac{(-1)^k a_k}{z^{k+1}} = \sum_{n=0}^{\infty} \frac{(-1)^n}{z(z+1)...(z+n)} \left\{ \sum_{k=0}^{n} s(n,k) \, a_k \right\}$$

replacing $a_k$ by $(-1)^k a_k$ and using (19).

The representations (13) follow from the identities

$$\sum_{k=0}^{n} \begin{bmatrix} n \\ k \end{bmatrix} = n!, \quad \sum_{k=0}^{n} \begin{bmatrix} n \\ k \end{bmatrix} k = n! H_n.$$

**Example 7**. In the same line, using the identity

$$\sum_{k=0}^{n} \begin{bmatrix} n \\ k \end{bmatrix} k^2 = n!\left(H_n + H_n^2 - H_n^{(2)}\right)$$

where $H_n^{(2)} = 1 + \dfrac{1}{2^2} + \ldots + \dfrac{1}{n^2}$, we find from (18) the representation

$$\sum_{k=0}^{\infty} \frac{k^2}{z^{k+1}} = \sum_{n=0}^{\infty} \frac{n!}{z(z+1)\ldots(z+n)} \left(H_n + H_n^2 - H_n^{(2)}\right).$$

Computing the series on the left hand side for $z > 1$ we come to the inverse factorial series

$$\frac{z+1}{(z-1)^3} = \sum_{n=0}^{\infty} \frac{n!}{z(z+1)\ldots(z+n)} \left(H_n + H_n^2 - H_n^{(2)}\right).$$

**Example 8**. In this example we use the Stirling numbers of the second kind $S(n,m)$ which can be defined by the generating function

(21) $\quad \displaystyle\sum_{n=0}^{\infty} S(n,m) x^n = \frac{x^m}{(1-x)(1-2x)\ldots(1-mx)}.$

For any two integers $0 \leq p \leq n$ the following identity is true (entry (A30) in [5])

$$\sum_{k=0}^{n} \begin{bmatrix} n \\ k \end{bmatrix} S(k+1, p+1) = \frac{n!}{p!} \binom{n}{p}.$$

According to Proposition 5 this implies the series identity

(22) $\quad \displaystyle\sum_{k=0}^{\infty} \frac{S(k+1, p+1)}{z^{k+1}} = \frac{1}{p!} \sum_{n=0}^{\infty} \frac{n!}{z(z+1)\ldots(z+n)} \binom{n}{p}.$

In view of (21) with $x = 1/z$ we find the representation

(23) $\quad \dfrac{1}{(z-1)(z-2)\ldots(z-p-1)} = \dfrac{1}{p!} \displaystyle\sum_{n=0}^{\infty} \dfrac{n!}{z(z+1)\ldots(z+n)} \binom{n}{p}$

convergent for $\text{Re}(z) > p + 1$. For $p = 0$ this is equation (13) and for $p = 1$ this is (15).

**Example 9**. Let $d_n$ be the sequence of Cauchy numbers of the second kind defined by the generating function

$$\frac{-x}{(1-x)\ln(1-x)} = \sum_{n=0}^{\infty} d_n \frac{x^n}{n!}$$

(see [9, p. 293]). The numbers $d_n$ for $n = 0, 1, \ldots$, satisfy the identity (entry (A42) in [5])

$$\sum_{k=0}^{n}\begin{bmatrix}n\\k\end{bmatrix}\frac{1}{k+1}=(-1)^{n}d_{n}$$

and from Proposition 5 it follows that

$$\sum_{k=0}^{\infty}\frac{1}{(k+1)z^{k+1}}=\sum_{n=0}^{\infty}\frac{(-1)^{n}d_{n}}{z(z+1)...(z+n)},\quad z>1.$$

The series on the left is easy to recognize and so we have the representation

(24) $$-\ln\left(1-\frac{1}{z}\right)=\sum_{n=0}^{\infty}\frac{(-1)^{n}d_{n}}{z(z+1)...(z+n)},\quad z>1.$$

It is interesting to compare this to the representation

(25) $$\ln\left(1+\frac{1}{z}\right)=\sum_{n=0}^{\infty}\frac{(-1)^{n}c_{n}}{z(z+1)...(z+n)},\quad z>0$$

where $c_n$ are the Cauchy numbers of the first kind defined by

$$\frac{x}{\ln(1+x)}=\sum_{n=0}^{\infty}c_{n}\frac{x^{n}}{n!}$$

(see [6, 9], and [30, p. 144]). A similar representation can be found on p. 244 in [20].

## 5. Asymptotic expansions

We can use Proposition 6 to obtain asymptotic series on the powers of $z^{-n-1}$ for functions with inverse factorial representations. We will use Proposition 5 to give proofs for the asymptotic representation of Nielsen's beta function $\beta(z)$ from Example 4, the incomplete gamma function $\gamma(z,x)$, and also $\psi'(z)$. We can use either equation (18) or equation (20).

Setting

$$\beta(z)=\sum_{n=0}^{\infty}\frac{n!}{2^{n+1}}\frac{1}{z(z+1)...(z+n)}=\sum_{n=0}^{\infty}\frac{1}{z(z+1)...(z+n)}\left\{\sum_{k=0}^{n}\begin{bmatrix}n\\k\end{bmatrix}a_{k}\right\}$$

we solve for $a_k$ from the equation

$$\sum_{k=0}^{n}\begin{bmatrix}n\\k\end{bmatrix}a_{k}=\frac{n!}{2^{n+1}}.$$

Using the formulas for the Stirling transform of sequences (see the Appendix in [5]) we have

$$(-1)^n a_n = \sum_{k=0}^{n} S(n,k) k! \frac{(-1)^k}{2^{k+1}} = \frac{1}{2} \sum_{k=0}^{n} S(n,k) k! \left(-\frac{1}{2}\right)^k.$$

At this point we involve the geometric polynomials

(26) $$\omega_n(x) = \sum_{k=0}^{n} S(n,k) k! x^k$$

introduced and studied in [7]. It is known that

$$\omega_n\left(-\frac{1}{2}\right) = \frac{2}{n+1}(1-2^{n+1}) B_{n+1} = E_n(0)$$

where $B_k$ are the Bernoulli numbers and $E_k(x)$ are the Euler polynomials. For the first equality see the solution of problem 6.76 on p. 559 in [10]. For the second equality see [3]. This way we have

$$a_n = \frac{(-1)^n}{n+1}(1-2^{n+1}) B_{n+1} = \frac{(-1)^n}{2} E_n(0)$$

and we come to following proposition.

**Proposition 7**. *Nielsen's beta function has the asymptotic series*

(27) $$\beta(z) = \sum_{n=0}^{\infty} \frac{(-1)^n (1-2^{n+1}) B_{n+1}}{n+1} \frac{1}{z^{n+1}} = \frac{1}{2} \sum_{n=0}^{\infty} (-1)^n E_n(0) \frac{1}{z^{n+1}}.$$

(Cf. Example 1 on p. 145 in [2]). In the same way we can use the representation of the incomplete gamma function $\gamma(z, x)$ from Example 5. Solving for the coefficients $a_n$ from the equation

(28) $$\sum_{k=0}^{n} \begin{bmatrix} n \\ k \end{bmatrix} a_k = x^n$$

we find that

(29) $$a_n = (-1)^n \sum_{k=0}^{n} S(n,k)(-1)^k x^k.$$

This time we use the exponential polynomials

$$\varphi_n(x) = \sum_{k=0}^{n} S(n,k) x^k$$

(see [3, 7]). These polynomials appeared in the works of Ramanujan [2], E.T. Bell, J. Touchard, Gian-Carlo Rota, and many others. As mentioned in [3], these polynomials were used as early as 1843 by the prominent German mathematician Johann A. Grunert.

Equation (29) says that $a_n = (-1)^n \varphi_n(-x)$, so that we have

$$\gamma(z, x) = x^z e^{-x} \sum_{n=0}^{\infty} (-1)^n \varphi_n(-x) \frac{1}{z^{n+1}}$$

in accordance with Proposition 5.1 in [7] proved there by a different method.

We notice that equations (28) and (29) in view of Proposition 5 lead to the interesting formula

(30) $$\sum_{n=0}^{\infty} \frac{x^n}{z(z+1)...(z+n)} = \sum_{k=0}^{\infty} (-1)^n \varphi_n(-x) \frac{1}{z^{k+1}}$$

which is Ramanujan's Example (c) on p.145 in [2] (see also p. 47 in [2]). In a similar manner, using the representation from Example 3

$$\psi'(z) = \sum_{n=0}^{\infty} \frac{n!}{n+1} \frac{1}{z(z+1)...(z+n)} = \sum_{n=0}^{\infty} \frac{(-1)^n n!}{n+1} \frac{(-1)^n}{z(z+1)...(z+n)}.$$

We compare this to (20) and solve for $a_k$ from the equation

$$\frac{(-1)^n n!}{n+1} = \sum_{k=0}^{n} s(n,k) a_k$$

to find

$$a_n = \sum_{k=0}^{n} S(n,k) \frac{(-1)^k k!}{k+1}$$

which are exactly the Bernoulli numbers [3], $a_n = B_n$. In view of (20) we find a new proof of the asymptotic expansion

(31) $$\psi'(z) = \sum_{k=0}^{\infty} \frac{(-1)^k B_k}{z^{k-1}}.$$

## 6. A binomial formula with inverse factorials

In this short section we prove a formula relating a binomial sum and an inverse factorial series with Stirling numbers of the first kind.

**Proposition 8**. *For any three integers $p \geq 1$, $k \geq 0$, $m \geq 0$ we have*

(32) $$\sum_{j=0}^{m} \binom{m}{j} \frac{(-1)^j}{(j+p)^{k+1}} = (p-1)! \sum_{n=0}^{\infty} \begin{bmatrix} n \\ k \end{bmatrix} \frac{1}{n!(n+m+1)(n+m+2)...(n+m+p)}.$$

*Proof*. We start with the (exponential) generating function for the unsigned Stirling numbers of the first kind

$$(-1)^k \ln^k(1-x) = k! \sum_{n=0}^{\infty} \begin{bmatrix} n \\ k \end{bmatrix} \frac{x^n}{n!}$$

where $|x| < 1$. We multiply both sides by $x^m(1-x)^{p-1}$ and integrate between $0$ and $1$ to get

(33) $\quad (-1)^k \int_0^1 x^m (1-x)^{p-1} \ln^k(1-x) dx = k! \sum_{n=0}^{\infty} \begin{bmatrix} n \\ k \end{bmatrix} \frac{1}{n!} \int_0^1 x^{n+m}(1-x)^{p-1} dx$.

In the integral on the left hand side we make the substitution $1-x = e^{-t}$

$$(-1)^k \int_0^1 x^m (1-x)^{p-1} \ln^k(1-x) dx = \int_0^{\infty} (1-e^{-t})^m e^{-pt+t} t^k e^{-t} dt$$

$$= \sum_{j=0}^m \binom{m}{j}(-1)^j \int_0^{\infty} t^k e^{-(j+p)t} dt = k! \sum_{j=0}^m \binom{m}{j}(-1)^j \frac{1}{(j+p)^{k+1}}.$$

For the integral on the right hand side we use the representation

$$\int_0^1 x^{n+m}(1-x)^{p-1} dx = \frac{(p-1)!}{(n+m+1)(n+m+2)...(n+m+p)}$$

which comes from Euler's Beta function

$$B(u,v) = \int_0^1 x^{u-1}(1-x)^{v-1} dx = \frac{\Gamma(u)\Gamma(v)}{\Gamma(u+v)}.$$

Equating both sides in (33) we come to the desired formula.

**Example 10**. In particular, for $p = 1$ in the above formula we have

(33) $\quad \sum_{j=0}^m \binom{m}{j} \frac{(-1)^j}{(j+1)^{k+1}} = \sum_{n=0}^{\infty} \begin{bmatrix} n \\ k \end{bmatrix} \frac{1}{n!(n+m+1)}.$

For $m = 0$, $m = 1$, and $m = 2$ this gives correspondingly

$$1 = \sum_{n=0}^{\infty} \begin{bmatrix} n \\ k \end{bmatrix} \frac{1}{n!(n+1)}$$

$$1 - \frac{1}{2^{k+1}} = \sum_{n=0}^{\infty} \begin{bmatrix} n \\ k \end{bmatrix} \frac{1}{n!(n+2)}$$

$$1 - \frac{1}{2^k} + \frac{1}{3^{k+1}} = \sum_{n=0}^{\infty} \begin{bmatrix} n \\ k \end{bmatrix} \frac{1}{n!(n+3)}$$

etc.